\documentclass[12pt,a4paper]{article}
\usepackage{amssymb}

\usepackage{graphicx,amssymb,amsfonts,epsfig,amsthm,a4,amsmath,url}
\usepackage[latin1]{inputenc}

\usepackage{enumitem}

\usepackage{amssymb}

\usepackage{graphicx,amssymb,amsfonts,epsfig,a4,amsmath,amsthm}
\usepackage[latin1]{inputenc}

\newtheorem{Thm}{Theorem}[section]
\newtheorem{Prop}[Thm]{Proposition}
\newtheorem{Lem}[Thm]{Lemma}
\newtheorem{Cor}[Thm]{Corollary}

\newtheorem{Rem}[Thm]{Remark}

\newcommand{\R}{\mathbf{R}}
\newcommand{\C}{\mathbf{C}}

\newcommand{\K}{\mathbf{K}}

\newcommand{\E}{\mathbf{E}}

\newcommand{\bprr}{\noindent \textbf{Proof} }
\newcommand{\epr}{~$\blacksquare$}

\title{The exceptional simple Lie group $F_{4(-20)}$, \\ after J. Tits}
\author{Alain VALETTE}

\begin{document}

\maketitle

\begin{abstract} This is a semi-survey paper, where we start by advertising Tits' synthetic construction from \cite{Tits}, of the hyperbolic plane $H^2(Cay)$ over the Cayley numbers $Cay$, and of its automorphism group which is the exceptional simple Lie group $G=F_{4(-20)}$. Let $G=KAN$ be the Iwasawa decomposition. Our contributions are:
\begin{itemize}
\item Writing down explicitly the action of $N$ on $H^2(Cay)$ in Tits'model, facing the lack of associativity of $Cay$.
\item If $MAN$ denotes the minimal parabolic subgroup of $G$, characterizing $M$ geometrically. 
\end{itemize}
 \end{abstract}


\section{Introduction}

Let $G$ be a simple Lie group with trivial center, with Iwasawa decomposition $G=KAN$. The {\it real rank} of $G$, i.e. the dimension of $A$, can be viewed as a measure of the complexity of $G$, with the real rank of $G$ being 0 if and only if $G$ is compact. It is known that in higher real rank, i.e. $\dim A\geq 2$, the group $G$ satisfies various forms of rigidity (e,g. property (T), arithmeticity of lattices, Margulis' super-rigidity: see \cite{Zim} for all that). The situation in real rank one is more subtle, with some forms of rigidity being satisfied and some other not, the situation depending very much on the group under consideration. Recall that, from Cartan's classification (see e.g. \cite{Tab}), simple real Lie groups with real rank 1 fall into 3 infinite series and one exceptional group:

\begin{itemize}
\item Adjoint group of $SO_0(n,1),\;(n\geq 2)$.
\item Adjoint group of $SU(n,1),\;(n\geq 2)$.
\item Adjoint group of $Sp(n,1),\;(n\geq 2)$.
\item $F_{4(-20)}$.
\end{itemize}

The present paper deals with the exceptional group $F_{4(-20)}$. Several papers in harmonic analysis on real rank 1 simple Lie groups deal exclusively with the classical cases, i.e. the three infinite series above: see \cite{Lip} for a sample. A notable exception is Takahashi's paper \cite{Taka}, which was inspiring for us. It is therefore interesting to have explicit constructions of $F_{4(-20)}$, both algebraic and geometric.

Algebraically, $F_{4(-20)}$ was constructed by Freudenthal \cite{Freud} as the automorphism group of the real Albert algebra $J_{1,2}$, which is one out of the three 27-dimensional exceptional Jordan algebras over the reals (see also chapter 19 in \cite{Mostow}, or more recently \cite{Allcock}). Geometrically, it was constructed by Tits \cite{Tits} as the automorphism group of the hyperbolic plane over the octonions $Cay$. This provides a unified construction of all simple Lie groups of real rank 1 and with trivial centre: let $\K$ be one of the four finite-dimensional division algebras over the reals, namely $\R,\C,\mathbf{H}$ (the Hamilton quaternions), or $Cay$ (the Cayley numbers). For $n\geq 2$, let $H^n(\K)$ be $n$-dimensional hyperbolic space over $\K$ (with $n=2$ when $\K=Cay$). Then any real rank 1 simple Lie group with trivial center is isomorphic to the (connected component of identity of the) automorphism group $Aut(H^n(\K))$ for some $n\geq 2$ and some $\K\in\{\R,\C,\mathbf{H},Cay\}$.

The present, mostly survey, paper aims to advertise Tits' synthetic construction of $H^2(\K)$ from \cite{Tits}. Tits first defines the projective plane $P^2(\K)$ as the affine plane $\K^2$ together with a line at infinity, which provides a very suitable coordinatization of $P^2(\K)$. The hyperbolic plane then appears as the set of interior points of a specific hyperbolic polarity, and is realized as the unit ball $\{(x_1,x_2): |x_1|^2+|x_2|^2<1\}$ in $\K^2$.

The paper is structured as follows: in section 2, we recall Cayley numbers and the triality principle; in section 3, we survey Tits' constructions and results from \cite{Tits}. Only the final section 4 has some novelty. Let $KAN$ be the Iwasawa decomposition of $F_{4(-20)}$; the maximal compact subgroup $K$ was determined by Tits in \cite{Tits}, the one-parameter subgroup $A$ appears in \cite{Taka}. The nilpotent subgroup $N$ was determined by Takahashi \cite{Taka} in Freudenthal's algebraic model, it is a central extension of the additive group of $Cay$ by the additive group $Im(Cay)$ of imaginary Cayley numbers. Our main contribution consists in writing down explicitly the action of $N$ on $H^2(Cay)$: doing so we pay the price for working over a non-associative division algebra. Indeed, here is for example the action of the element $u(y,0)$ in the centre $Z(N)$ on $(x_1,x_2)\in H^2(Cay)$ (see Corollary \ref{actionofZ(N)} below): we have $u(y,0)(x_1,x_2)=(x'_1, x'_2)$, where:
$$x'_1=((1-y)x_1 + y)(-yx_1 +1+y)^{-1} $$
and
$$x'_2=(x_2(1-\overline{x_1}))[(1-\overline{x_1})^{-1}(-yx_1+1+y)^{-1}]$$
While the formula for $x'_1$ is what is expected for an action by M\"obius transformations, the formula for $x'_2$ reveals the problems that arise with the lack of associativity.

Finally, let $MAN$ be the minimal parabolic subgroup of $F_{4(-20)}$: we characterize $M$ as the stabilizer in $K$ of the point $(1,0)$ on the boundary of $H^2(Cay)$.


\section{Octonions and triality}

Recall that $Cay$, the algebra of {\it Cayley numbers} or {\it octonions}, denotes the 8-dimensional non-associative division algebra over $\R$ with basis $e_0, e_1, e_2,..., e_7$
and relations
\begin{itemize}
\item $e_0.e_i=e_i.e_0=e_i$ for $i=0,1,...,7$.
\item $e_i^2=-e_0$ for $i>0$.
\item $e_i.e_j=-e_j.e_i$ for $i,j>0,\;i\neq j$.
\item $e_2.e_6 = e_3.e_4= e_5.e_7=e_1$, plus all relations deduced from this one by cyclic permutations of the indices $1,2,...,7$.
\end{itemize}

We will replace $e_0$ by 1 for simplicity and identify the real octonion $a.e_0$ with the real number $a$. Any octonion $x=a + \sum_{i=1}^7 b_ie_i$ has a {\it conjugate} $\overline{x}=a-\sum_{i=1}^7 b_ie_i$ and a {\it norm} $N(x)=x\overline{x}=\overline{x}x=a^2+ \sum_{i=1}^7 b_i^2$. The {\it modulus} of $x$ is $|x|=:\sqrt{N(x)}$. The norm satisfies $N(xy)=N(x)N(y)$ for every $x,y\in Cay$. 

The {\it imaginary part} of an octonion $x$ is $Im(x)=\frac{x-\overline{x}}{2}$: the octonion $x$ is said to be {it imaginary} if $x=Im(x)$; we denote by $Im(Cay)$ the set of imaginary octonions.

The algebra $Cay$ is a division algebra as the inverse of a non-zero octonion $x$ is $x^{-1}=\frac{\overline{x}}{N(x)}$. The algebra $Cay$ is not associative but {\it alternative}, namely for every $x,y\in Cay$:
\begin{itemize}
\item  $(xy)x=x(yx)$;
\item $x(xy)=x^2y$;
\item $(xy)y=xy^2$.
\end{itemize}
For all this, see either the nice paper \cite{Baez} or Chapter 1 of the book \cite{SpVe}.

The norm allows us to identify $Cay$ linearly isometrically with the 8-dimensional Euclidean space $\E^8$. So the rotation group $SO(8)$ acts linearly on $Cay$ and the have the {\it triality principle} (see Theorem 3.2.1 in \cite{SpVe}) that we recall:

\begin{Thm}\label{trial} For every $R\in SO(8)$, there exists $R',R''\in SO(8)$ such that, for every $x,y\in Cay$:
$$R(xy)=R'(x)R''(y)$$
Moreover $R',R''$ are determined up to sign by this relation.
\epr
\end{Thm}

We will need some special features of the triple $(R,R',R'')$ as above. 

\begin{Lem}\label{monlemme} Fix $R\in SO(8)$, let $R',R''\in SO(8)$ be associated with $R$ as in Theorem \ref{trial}.
\begin{enumerate}
\item We have $R''(1)=1$ if and only if $R=R'$.
\item Assume $R''(1)=1$. Then:
\begin{enumerate}
\item[a)] For every non-zero octonion $x$, we have $R''(x^{-1})=R''(x)^{-1}$.
\item[b)] $R''$ is a Jordan homomorphism, i.e 
$$R''(xy+yx)=R''(x)R''(y)+R''(y)R''(x)$$
for every $x,y\in Cay$.
\item[c)] For every $x\in Cay$, the map $R''$ is a homomorphism in restriction to the subalgebra generated by $1$ and $x$.
\end{enumerate}
\end{enumerate}
\end{Lem}

\bprr: \begin{enumerate}
\item If $R''(1)=1$, then for every $x\in Cay$ we get:
$$R(x)=R(x.1)=R'(x).R''(1)=R'(x).1=R'(x)$$
hence $R=R'$. Conversely, assuming $R=R'$, we get for every $x\in Cay$:
$$R(x)=R(x.1)=R'(x)R''(1)=R(x)R''(1),$$
which forces $R''(1)=1$.

\item \begin{enumerate}
\item[a)] As $R''(1)=1$, we have $R''(\overline{x})=\overline{R''(x)}$, and the result follows from the formula for the inverse of $x$.
\item[b)] Denoting by $\langle .|.\rangle$ the standard scalar product on $Cay$ associated with the norm $N$, we have by formula (1.8) in \cite{SpVe}:
$$xy+yx= \langle x|1\rangle y+ \langle y|1\rangle x - \langle x|y\rangle 1.$$
Applying $R''$ and using that it is a linear isometry, we get:
$$R''(xy+yx)=\langle R''(x)|R''(1)\rangle R''(y)+\langle R''(y)| R''(1)\rangle R''(x)- \langle R''(x)|R''(y)\rangle R''(1)$$
$$= \langle R''(x)|1\rangle R''(y)+\langle R''(y)|1\rangle R''(x) - \langle R''(x)|R''(y)\rangle 1 = R''(x)R''(y)+R''(y)R''(x),$$
where the second equality follows from the assumption $R''(1)=1$.
\item[c)] By formula (1.7) in \cite{SpVe}, the subalgebra generated by $1$ and $x$ is at most 2-dimensional, so it is commutative and the result follows from the preceding point.
\epr
\end{enumerate}
\end{enumerate}

\section{The hyperbolic plane over $Cay$}

In this section we follow rather closely Tits'paper \cite{Tits}.

\subsection{The projective plane and its collineation group}

Let $\K$ denote one of the following division algebras over the reals: $\R,\C,\mathbf{H}$(the Hamilton quaternions), or $Cay$. The {\it projective plane} $P^2(\K)$ is the disjoint union
$$P^2(\K)=(\K\times \K)\coprod \K\coprod{(\infty)}$$
where $\K\times \K$ is the set of points at finite distance and $\K\coprod{(\infty)}=\{(u):u\in \K\coprod\{\infty\}\}$ is the set of point at infinity. A {\it line} in $P^2(\K)$ is a subset of one of the following three types:
\begin{itemize}
\item $[u,v]= \{(x,y)\in \K\times \K: y=ux+v\}\cup\{(u)\}$;
\item $[u]=\{(x,y)\in \K\times \K: x=u\}\cup\{(\infty)\}$;
\item $[\infty]=\K\coprod{(\infty)}$.
\end{itemize}
Lines of the first type are oblique, lines of the second type are vertical, and $[\infty]$ is the line at infinity. It is readily seen that the incidence axioms of projective geometry are satisfied: two lines intersect in a unique point, two distinct points determine a unique line through them. 

A {\it collineation} of $P^2(\K)$ is a bijection of $P^2(\K)$ onto itself, that preserves the set of lines. 

\begin{Thm} (Th\'eor\`eme 4.9 in \cite{Tits}) The group of collineations of $P^2(Cay)$ is isomorphic to the real form $E_{6(-26)}$ of the exceptional Lie group $E_6$: it is therefore a simple Lie group of dimension 78, which is simply connected with trivial centre. 
\epr
\end{Thm}

\subsection{The hyperbolic plane and its automorphism group}

A {\it polarity} of $P^2(\K)$ is a bijection $\Pi$ between the set of points and the set of lines of $P^2(\K)$, satisfying two conditions:
\begin{itemize}
\item if $a$ runs through a line $L$, then $\Pi(a)$ runs in the sheaf of lines through a given point;
\item for any two points $a,b\in P^2(\K)$, we have $a\in \Pi(b)\Leftrightarrow b\in\Pi(a)$.
\end{itemize}

We will need just one example of polarity, given on points by:
$$\Pi(x,y)=[-\overline{xy^{-1}},\overline{y^{-1}}]\;\;\mbox{if}\,y\neq 0;$$
$$\Pi(x,0)=[\overline{x^{-1}}]\;\;\mbox{if}\;x\neq 0;$$
$$\Pi(0,0)=[\infty].$$

This polarity is {\it hyperbolic}, meaning that the set of points $P\in P^2(\K)$ such that $P\in\Pi(P)$ is non-empty. Actually this set is the sphere of equation $|x|^2+|y|^2=1$ in $\K\times \K$. The {\it hyperbolic plane} $H^2(\K)$ is then the ball $|x|^2+|y|^2<1$ in $\K\times \K$. Lines in $H^2(\K)$ are traces on $H^2(\K)$ of lines in $P^2(\K)$ that intersect it (``Klein model''). An automorphism of $H^2(\K)$ is a collineation of $P^2(\K)$ that commutes with $\Pi$.

\begin{Thm}\label{TitsF_4}  Let $G$ be the automorphism group of $H^2(Cay)$.
\begin{enumerate}
\item (6.5 in \cite{Tits}) The stabilizer $K$ of $(0,0)$ in $G$ is isomorphic to $Spin(9)$ (the universal cover of $SO(9)$); it acts transitively on the sphere $|x|^2+|y|^2=1$ and it is a maximal compact subgroup of $G$.
\item (Th\'eor\`eme 6.8 in \cite{Tits}) $G$ is isomorphic to the real form $F_{4(-20)}$ of the exceptional Lie group $F_4$: it is therefore a simple Lie group of dimension 52, which is simply connected with trivial centre. \epr
\end{enumerate}
\end{Thm}

Let us single out the following for future use:

\begin{Lem}\label{LemTits} (see (4.6) in \cite{Tits}) An element $k\in K$ fixes the points at infinity $(0)$ and $(\infty)$ if and only if there exists $R\in SO(8)$ such that (in the notation of the triality principle Theorem \ref{trial}):
$$k(x,y)=(R''(x), R(y))$$
for every $(x,y)\in Cay\times Cay$. \epr
\end{Lem}

\section{The Iwasawa decomposition and the minimal parabolic subgroup of $G$}

Let $G$ and $K$ be as in Theorem \ref{TitsF_4}. Let $\mathfrak{g,k}$ denote the Lie algebras of $G,\,K$ respectively. Using the algebraic realization of $G$ as the automorphism group of the Albert algebra $J_{1,2}$, Takahashi (\cite{Taka}, formula (10) in section 3) determines the Cartan decomposition $\mathfrak{g}=\mathfrak{k}\oplus\mathfrak{p}$, chooses a natural maximal abelian subalgebra $\mathfrak{a}\subset\mathfrak{p}$, where $\mathfrak{a}=\R\cdot A$ (as $\dim\mathfrak{a}=1$) and identifies restricted root subspaces 
$$\mathfrak{g}_\alpha=\{Z(z):z\in Cay\}$$
and
$$\mathfrak{g}_{2\alpha}=\{Y(y): y\in Im(Cay)\}.$$
Setting $\mathfrak{n}=\mathfrak{g}_\alpha\oplus\mathfrak{g}_{2\alpha}$ he establishes the Iwasawa decomposition at the level of $\mathfrak{g}$ (\cite{Taka}, (13) in section 3):
$$\mathfrak{g}=\mathfrak{k}\oplus\mathfrak{a}\oplus\mathfrak{n}.$$

We then set $a_t=\exp(tA)$ and $A=\{a_t:t\in\R\}$, the split maximal torus $A$ in the Iwasawa decomposition of $G$. Similarly we set
$$u(y,z)=\exp(Y(y))\exp(Z(z))=\exp(Z(z))\exp(Y(y))=\exp(Y(y)+Z(z))$$
and $N=\{u(y,z): y\in Im(Cay), z\in Cay\}$, the nilpotent part in the Iwasawa decomposition of $G$; so $N$ is a 15-dimensional, 2-step nilpotent group with centre $Z(N)=\{u(y,0):y\in Im(Cay)\}$.

In order to write down the action of $A$ and $N$ on $H^2(Cay)$, Takahashi defines in section 5 of \cite{Taka} a $G$-equivariant embedding of $H^2(Cay)$ in $J_{1,2}$, which allows him to write down the action of $A$ in full generality, but the action of $N$ only at the origin $(0,0)\in H^2(Cay)$:


\begin{Thm} (formulae (14) and (20) in section 5 of \cite{Taka}) 
\begin{enumerate} 
\item For $(x_1,x_2)\in H^2(Cay),t\in\R$, we have 
\begin{equation}
a_t(x_1,x_2)= ((\cosh(t)x_1 + \sinh(t))(\sinh(t)x_1+\cosh(t))^{-1}, x_2(\sinh(t)x_1+\cosh(t))^{-1}).
\end{equation}
\item For $t\in\R,y\in Im(Cay),z\in Cay$, we have $a_tu(y,z)(0,0) = (x_1,x_2)$, where:
\begin{equation}
x_1= (\sinh(t) + e^t(\frac{|z|^2}{2}+y))(\cosh(t) + e^t(\frac{|z|^2}{2}+y))^{-1}
\end{equation}
\begin{equation}
x_2=z(\cosh(t) + e^t(\frac{|z|^2}{2}+y))^{-1}
\end{equation}\epr
\end{enumerate}
\end{Thm}

Observe that formula (1) is the expected formula for a hyperbolic M\" obius transformation. Its validity is due to the fact that coefficients of $a_t$ are real, hence belong to the centre of $Cay$. Our aim now is to jack up formulae (2) and (3) to get the action of $N$ on any point $(x_1,x_2)$ in $H^2(Cay)$.

\begin{Thm}\label{actionofN} For $(x_1,x_2)\in H^2(Cay), y\in Im(Cay), z\in Cay$, we have $u(y,z)(x_1,x_2)=(x'_1,x'_2)$ where:
\begin{equation}
x'_1= [((1-\frac{|z|^2}{2}-y)x_1+ (\overline{z}(x_2(1-\overline{x_1})))(1-\overline{x_1})^{-1}+ \frac{|z|^2}{2}+y)(1-\overline{x_1})].D^{-1},
\end{equation}
\begin{equation}
x'_2= [(-zx_1+x_2+z)(1-\overline{x_1})].D^{-1},
\end{equation}
and
\begin{equation}
D=[-(\frac{|z|^2}{2}+y)x_1 + (\overline{z}(x_2(1-\overline{x_1})))(1-\overline{x_1})^{-1} +1+\frac{|z|^2}{2}+y](1-\overline{x_1})
\end{equation}
\end{Thm}

\bprr We will need two rules for calculation in the group $AN$. First, by definition of $u(y,z)$ we have:
\begin{equation}
a_t u(y,z)a_{-t}=u(e^{2t}y,e^t{z}).
\end{equation}
Second, for $y_1,y_2\in Im(Cay),\,z_1,z_2\in Cay$, by lemma 4.5 in \cite{Nishio}:
\begin{equation}
u(y_1,z_1)u(y_2,z_2)=u(y_1+y_2+Im(\overline{z_1}z_2),z_1+z_2).
\end{equation}
 
Takahashi proves (see lemme 2 in section 5 of \cite{Taka}) that the group $AN$ acts simply transitively on $H^2(Cay)$. So, for $(x_1,x_2)\in H^2(Cay)$, there exist uniquely determined $t_0\in\R, y_0\in Im(Cay), z_0\in Cay$ such that
\begin{equation}
(x_1,x_2)=a_{t_0}u(y_0,z_0)(0,0),
\end{equation}
where, setting $r=(1-|x_1|^2-|x_2|^2)^{1/2}$, the objects $t_0, y_0, z_0$ are explicitly given (see formula (21) in section 5 of \cite{Taka}):
\begin{equation}
e^{-t_0}=\frac{|1-x_1|}{r}
\end{equation}
\begin{equation}
y_0=\frac{x_1-\overline{x_1}}{2r^2}
\end{equation}
\begin{equation}
z_0=\frac{x_2(1-\overline{x_1})}{|1-x_1|r}
\end{equation}

\medskip
Now we have by formulae (9), (7) and (8):
$$(x'_1,x'_2)=u(y,z)(x_1,x_2)=u(y,z)a_{t_0}u(y_0,z_0)(0,0)=a_{t_0}u(e^{-2t_0}y,e^{-t_0}z)u(y_0,z_0)(0,0)$$
$$=a_{t_0}u(e^{-2t_0}y+y_0+ e^{-t_0}Im(\overline{z}z_0),e^{-t_0}z+z_0)(0,0)$$
So by formulae (2) and (3):
\begin{equation}
x'_1= (\sinh(t_0)+E)(\cosh(t_0)+E)^{-1}
\end{equation}
\begin{equation}
x'_2=(e^{-t_0}z+z_0)(\cosh(t_0)+E)^{-1}
\end{equation}
where
\begin{equation}
E=e^{t_0}(\frac{|e^{-t_0}z+z_0|^2}{2}+ e^{-2t_0}y+ y_0+e^{-t_0}Im(\overline{z}z_0))
\end{equation}
It remains to express $x'_1,x'_2$ in terms of $x_1,x_2$. In the following computations, we will frequently use the following weak form of associativity in $Cay$: for every $a,b\in Cay$, we have, by lemma 1.3.3 in \cite{SpVe}:
$$(ab)\overline{b}=N(b)a.$$
We start by checking formula (5). By formulae (10) and (12):
$$e^{-t_0}z+z_0=\frac{|1-x_1|z}{r}+\frac{x_2(1-\overline{x_1})}{r|1-x_1|}=\frac{|1-x_1|^2z+x_2(1-\overline{x_1})}{r|1-x_1|}$$
$$=\frac{(z(1-x_1)+x_2)(1-\overline{x_1})}{r|1-x_1|}$$
Hence
\begin{equation}
e^{-t_0}z+z_0=\frac{(-zx_1+x_2+z)(1-\overline{x_1})}{r|1-x_1|}
\end{equation}
On the other hand by (15):
$$E=e^{t_0}[\frac{e^{-2t_0}|z|^2 + e^{-t_0}\overline{z}z_0+e^{-t_0}\overline{z_0}z+ |z_0|^2}{2} + e^{-2t_0}y+ y_0+e^{-t_0}(\frac{\overline{z}z_0-\overline{z_0}z}{2})]$$
$$= \frac{e^{-t_0}|z|^2+2\overline{z}z_0+e^{t_0}|z_0|^2+2e^{-t_0}y+2e^{t_0}y_0}{2};$$
hence by (10), (11), (12):
\begin{equation}
E=\frac{|1-x_1|^2|z|^2+2\overline{z}(x_2(1-\overline{x_1}))+|x_2|^2+2y|1-x_1|^2 + x_1-\overline{x_1}}{2r|1-x_1|}.
\end{equation}

Now by (10) we have, since $r^2=1-|x_1]^2-|x_2|^2$:
$$\cosh(t_0)=\frac{1}{2}(\frac{r}{|1-x_1|}+\frac{|1-x_1|}{r})=\frac{r^2+|1-x_1|^2}{2r|1-x_1|}=\frac{2-x_1-\overline{x_1}-|x_2|^2}{2r|1-x_1|}.$$
Hence by (17):
$$\cosh(t_0)+E= \frac{|1-x_1|^2|z|^2+2\overline{z}(x_2(1-\overline{x_1}))+2y|1-x_1|^2 + 2-2\overline{x_1}}{2r|1-x_1|}$$
Factoring out $2(1-\overline{x_1)}$ on the right in the numerator, we get:
$$\cosh(t_0)+E= \frac{[\frac{|z|^2}{2}(1-x_1)+(\overline{z}(x_2(1-\overline{x_1})))(1-\overline{x_1})^{-1}+y(1-x_1) + 1](1-\overline{x_1})}{r|1-x_1|}$$
hence
$$\cosh(t_0)+E= \frac{[-(\frac{|z|^2}{2}+y)x_1+(\overline{z}(x_2(1-\overline{x_1})))(1-\overline{x_1})^{-1}+1+\frac{|z|^2}{2}+y](1-\overline{x_1})}{r|1-x_1|}.$$
Comparing with (6), we get
\begin{equation}
\cosh(t_0)+E= \frac{D}{r|1-x_1|}.
\end{equation}
In view of (14), formula (5) follows immediately by combining (16) and (18).

\medskip
We now turn to verifying formula (4). Again from (10), we have:
$$\sinh(t_0)=\frac{1}{2}(\frac{r}{|1-x_1|}-\frac{|1-x_1|}{r})=\frac{r^2-|1-x_1|^2}{2r|1-x_1|}=\frac{x_1+\overline{x_1}-2|x_1|^2-|x_2|^2}{2r|1-x_1|},$$
hence by (17):
$$\sinh(t_0)+E= \frac{ |1-x_1|^2|z|^2 + 2\overline{z}(x_2(1-\overline{x_1})) +2y|1-x_1|^2 + 2x_1-2|x_1|^2}{2r|1-x_1|}.$$
Factoring out $2(1-\overline{x_1)}$ on the right in the numerator, we get:
$$\sinh(t_0)+E= \frac{[ \frac{|z|^2}{2}(1-x_1) + (\overline{z}(x_2(1-\overline{x_1})))(1-\overline{x_1})^{-1} +y(1-x_1) + x_1](1-\overline{x_1})}{r|1-x_1|}$$
hence
\begin{equation}
\sinh(t_0)+E= \frac{[ (1-\frac{|z|^2}{2}-y)x_1 + (\overline{z}(x_2(1-\overline{x_1})))(1-\overline{x_1})^{-1} +\frac{|z|^2}{2}+y](1-\overline{x_1})}{r|1-x_1|}
\end{equation}
In view of (13), formula (4) follows immediately by combining (18) and (19).
\epr

\medskip

In restriction to the center $Z(N)$, we get: 

\begin{Cor}\label{actionofZ(N)} For $(x_1,x_2)\in H^2(Cay), y\in Im(Cay)$, we have $u(y,0)(x_1,x_2)=(x'_1,x'_2)$ where:
\begin{equation}\label{expected}
x'_1=((1-y)x_1 + y)(-yx_1 +1+y)^{-1}
\end{equation}
\begin{equation}\label{unexpected}
x'_2=(x_2(1-\overline{x_1}))[(1-\overline{x_1})^{-1}(-yx_1+1+y)^{-1}]
\end{equation}
\end{Cor}

\bprr Formula (\ref{unexpected}) follows immediately from (5) and (6) by setting $z=0$. To obtain formula (\ref{expected}) in the same way, we must justify why we may use associativity to cancel out $1-\overline{x_1}$. The reason is that subalgebras of $Cay$ generated by 2 elements are associative (see Theorem 1.4.3 in \cite{SpVe}), and $\overline{x_1}$ belongs to the subalgebra generated by $x_1$ since $\overline{x_1}=|x_1|^2x_1^{-1}$. \epr

\begin{Rem}
Let us observe that, if we replace the Cayley numbers by the Hamilton quaternions $\mathbf{H}$, formulae (4) and (5) simplify drastically thanks to associativity:
$$x'_1= ((1-\frac{|z|^2}{2}-y)x_1 +\overline{z}x_2+ \frac{|z|^2}{2}+y)(-(\frac{|z|^2}{2}+y)x_1+\overline{z}x_2 + 1 +\frac{|z|^2}{2} + y)^{-1},$$
$$x'_2= (-zx_1+x_2+z)(-(\frac{|z|^2}{2}+y)x_1+\overline{z}x_2 + 1 +\frac{|z|^2}{2} + y)^{-1}$$
for $y\in Im(\mathbf{H}), z\in \mathbf{H}$.
This is nothing but the classical action of the nilpotent part of $Sp(2,1)$ on $H^2(\mathbf{H})$ by M\"obius tranformations (see e.g. page 21 in \cite{Lip}).
\end{Rem}

\medskip
We now turn to the {\it minimal parabolic subgroup} of $G$ i.e. the subgroup $P=MAN$, where $M$ denotes the commutant of $A$ in $K$. 

\begin{Prop}\label{montruc} We have $M=Stab_K(1,0)$, the stabilizer of $(1,0)$ in $K$. 
\end{Prop}

\bprr : Observe that $(1,0)$ is the unique attracting point of $A$, i.e. 
$$(1,0)=\lim_{t\rightarrow+\infty}a_t(x,y)$$ 
for every $(x,y)\in H^2(Cay)$. So for $m\in M$ we have
$$m(1,0)=\lim_{t\rightarrow+\infty}ma_t(x,y)=\lim_{t\rightarrow+\infty}a_tm(x,y)=(1,0),$$
so that $M\subset Stab_K(1,0)$. To prove the converse inclusion, take $k\in Stab_K(1,0)$. As $k$ stabilizes the line $[0,0]$ through $(0,0)$ and $(1,0)$, it fixes its point $(0)$ at infinity. It also fixes $\Pi[0,0]=(\infty)$. So lemma \ref{LemTits} applies, and there exists $R\in SO(8)$ such that $k(x,y)=(R''(x),R(y))$. Since $k(1,0)=(1,0)$, we have $R''(1)=1$. We now check that $k$ commutes with $a_t$: indeed, for every $(x,y)\in H^2(Cay)$:
$$a_tk(x,y)=((\cosh(t)R''(x)+\sinh(t))(\sinh(t)R''(x)+\cosh(t))^{-1}, R(y)(\sinh(t)R''(x)+\cosh(t))^{-1})$$
$$=(R''((\cosh(t)x + \sinh(t))(\sinh(t)x+\cosh(t))^{-1}),R'(y)R''((\sinh(t)x+\cosh(t))^{-1})\;\mbox{by lemma \ref{monlemme};}$$
$$=(R''((\cosh(t)x + \sinh(t))(\sinh(t)x+\cosh(t))^{-1}), R(y(\sinh(t)x+\cosh(t))^{-1}))\;\mbox{by Theorem \ref{trial};}$$
$$=ka_t(x,y).$$
This concludes the proof. \epr

\medskip
It was proved by Borel \cite{Bor} and reproved by Tits (see (6.9) in \cite{Tits}) that $Stab_K(1,0)$ is isomorphic to $Spin(7)$, the universal cover of $SO(7)$. 





\noindent
Author's address:\\
Institut de Math\'ematiques\\
Universit\'e de Neuch\^atel\\
Unimail\\
11 Rue Emile Argand\\
CH-2000 Neuch\^atel - SWITZERLAND

\medskip
\noindent
alain.valette@unine.ch

\end{document}